\documentclass[12pt]{article}
\usepackage{graphicx}
\usepackage{amssymb}
\usepackage{epstopdf}
\title{Corrig\' e}
\author{}
   
\usepackage{graphicx}
\usepackage{amssymb}
\usepackage{epstopdf}
\title {A family of two dimensionally determined ergodic processes }
\date{}
\author {Doureid Hamdan} 
\title
  {Ergodic quasi-exchangeable stationary processes  are isomorphic to  Bernoulli processes}
\begin{document}
\maketitle

{\abstract{\textwidth=4,5 in} A discrete time process, with law $\mu$, is quasi-exchangeable if for any  finite permutation $\sigma$  of time indices,  the law $\mu_\sigma$ of the resulting process is equivalent to $\mu$. For a quasi-exchangeable stationary process we prove mainly (1) that if the process is ergodic then it is isomorphic to a Bernoulli process and (2)  that if the family of all  Radon-Nikodym derivatives $\{{d\mu_\sigma\over d\mu}\}$ is uniformly integrable then the process is a mixture of  Bernoulli processes, which generalizes De Finetti's Theorem. We give application of (1)  to some determinantal processes. }


\section{Introduction, Notation}
According to the generalization, by Hewitt and Savage, of De Finetti's classical Theorem, 
 an exchangeable sequence $(X_n)_{n\ge 1}$ of  random variables, with values in a presentable space,  is a mixture of   sequences of independent identically distributed random variables. ( The family of presentable spaces include the polish or locally compact spaces).\\
 In the present paper we consider, more generally,  the class of   stationary quasi-exchangeable (  quasi symmetric, or quasi-invariant are also used) sequences of random variables (  Recall that  a sequence $(X_n)_{n\ge 1}$ of random variables is exchangeable [ respectively quasi-exchangeable] if its law $P$ is equal [ respectively equivalent ] to the law $P_\sigma$ of the sequence $(X_{\sigma(n)})_{n\ge 1}$, for  any permutation $\sigma$  of finitely many coordinates. \\  For a stationary quasi-exchangeable sequence  $(X_n)_{n\ge 1},$ 
 we prove mainly two facts:\\
(1) In general, we prove, Theorem 1, that  if the  dynamical system generated by $(X_n)_{n\ge 1}$ is ergodic, then it  is  isomorphic to a Bernoulli system.\\  
    The result (1) implies, Theorem 3,   that the discrete time  stationary   quasi-invariant determinantal   processes are  isomorphic to Bernoulli processes. These determinantal processes contain  the  discrete time  translation-invariant determinantal processes  in the family  considered by Bufetov in [3] (  for example,  the discrete sine process of Borodin Okounkov and Olshanski).\\
    Also, we give( Corollary 2, Remark 3) a simple proof of De Finetti's Theorem, when the state space $K$ is compact.\\
(2) 
Under the additional hypothesis that  the family  $\hskip 0,1 cm  \{ {dP_\sigma \over dP }:\sigma \textit {\hskip 0,1 cm is a finite permutation of } \hskip 0,1 cm \mathbb N\}\hskip 0,1 cm $  of Radon-Nikodym derivatives,   be uniformly integrable, we prove( Theorem 2)  that it is a mixture of i.i.d.  sequences. This generalizes the De Finetti's Theorem. \\ 
We establish, first, some definitions and notations.\\
{\bf{Definition 1}}\\
\textit{A sequence  $\hskip 0,1 cm (X_n)_{n\ge 1}\hskip 0,1 cm $ of random variables  is exchangeable if   the law $P_\sigma$ of the process $(X_{\sigma(n)})_{n\ge 1}$ is equal  the law $P$ of the process $(X_n)_{n\ge 1}$, for every permutation $\sigma$ belonging to the group $H_1$ of all permutations of the set $\mathbb N$ of natural numbers, which leave fixed all but a finite number of integers.}\\
  
Suppose that for any   $n\ge 1$, $X_n$ takes values in the measurable space $(K,{\cal F})$. 
  Then  the sequence  $(X_n)_{n\ge 1}$ is exchangeable, if and only if
  \begin{eqnarray}
P(X_1\in A_1,...,X_n\in A_n)=P(X_{\tau(1)}\in A_1,...,X_{\tau(n)}\in A_n),
\end{eqnarray} 
holds for all $\hskip 0,1 cm n\ge 1$, $A_1,...,A_n\in {\cal F}\hskip 0,1 cm $ and any permutation $\hskip 0,1 cm \tau\hskip 0,1 cm $ of $\hskip 0,1 cm \{1,...,n\}$,
or  equivalently, if and only if for any permutation $\hskip 0,1 cm \sigma\hskip 0,1 cm $ of $\hskip 0,1 cm \{1,...,n\}$ ( $\tau=\sigma^{-1}$),
\begin{eqnarray}
P(X_1\in A_1,...,X_n\in A_n)=P(X_1\in A_{\sigma(1)},...,X_n\in A_{\sigma(n)}).
\end{eqnarray}
In the particular case where  $\hskip 0,1 cm \tau\hskip 0,1 cm$ is defined by \begin{eqnarray*}\hskip 0,1 cm \tau(k)=k+1, \hskip 0,1 cm \textit{for} \hskip 0,1 cm \hskip 0,1 cm 1\le k\le n-1, \hskip 0,2 cm \textit{and} \hskip 0,3 cm \tau(n)=1,\hskip 0,1 cm\end{eqnarray*} and when $\hskip 0,1 cm A_n=K,\hskip 0,1 cm$  equality (1) reads
\begin{eqnarray}
P(X_1\in A_1,...,X_{n-1}\in A_{n-1})=P(X_2\in A_1,...,X_n\in A_{n-1})\hskip 0,7 cm (\textit{shift-inv}).
\end{eqnarray}
which  proves that the  law $P$  of an exchangeable sequence is  invariant by the unilateral shift on $(K^\mathbb N,{\cal F}^{\otimes \mathbb N} )$.\\
  
 Let  
\begin{eqnarray*}
\Omega:=K^\mathbb Z,\hskip 0,1 cm \textit{ endowed with the product sigma algebra} \hskip 0,3 cm {\cal B}:={\cal F}^{\otimes \mathbb Z},\hskip 0,2 cm  \textit{  and}\\\hskip 0,2 cm       X'_n(\omega)=\omega_n,\hskip 0,2 cm  \textit{ for all integer }\hskip 0,05 cm n\in \mathbb Z\hskip 0,3 cm \textit{and all}\hskip 0,4 cm \omega\in \Omega. \end{eqnarray*}
Suppose that $(X_n)_{n\ge 1}$, with law $P$, is exchangeable.  Define the measure $\mu$ on $(\Omega,{\cal B})$, which extends $P$,  by setting  for all $l,k\ge 0$, and $\hskip 0,1 cm A_{-l},...,A_k\in {\cal F},\hskip 0,1 cm$      
\begin{eqnarray*}
\mu(X'_{-l}\in A_{-l},...,X'_k\in A_k):=P(X_1\in A_{-l},...,X_{k+l+1}\in A_k).
\end{eqnarray*}
Then, by (3), for all $p\ge 1$,
\begin{eqnarray}
\mu(X'_{-l}\in A_{-l},...,X'_k\in A_k)=P(X_p\in A_{-l},...,X_{k+l+p}\in A_k).\hskip 0,7 cm (\textit{shift-inv})'.
\end{eqnarray}
so that  $\mu$ extends to a probability measure  on $\Omega$, which   is also invariant by  the shift transformation $\hskip 0,1 cm S:\hskip 0,1 cm (S\omega)_n=\omega_{n+1},\hskip 0,1 cm n\in  \mathbb Z,\hskip 0,1 cm \omega \in \Omega. $ \\
Let   $G$ be the group of all permutations of $\hskip 0,1 cm \mathbb Z\hskip 0,1 cm $  and
 $\hskip 0,1 cm H\subset G,$ be  the subgroup of all permutations  with finite support:  \begin{eqnarray*}
\sigma\in H\iff \sigma\in G, \hskip 0,2 cm \textit{and} \hskip 0,3 cm  \exists N, \hskip 0,2 cm  \sigma(n)=n, \forall n,  \mid n\mid \ge N.\end{eqnarray*}
For any $\tau\in G, \hskip 0,2 cm$ let  the transformation $T_\tau:\Omega\rightarrow \Omega$, be defined for all  $ \hskip 0,2 cm \omega\in \Omega, \hskip 0,2 cm$  by 
 \begin{eqnarray*}(T_\tau(\omega))_n=\omega_{\tau(n)}, \hskip 0,2 cm \forall  n\in \mathbb Z.\end{eqnarray*} 
 Then $T_\tau$ is ${\cal B}$-measurable, and, when $K$ is a topological space,  $T_\tau$ is continuous for the product topology on  $\Omega$. Also     for all $\sigma$ and $\tau$ in $H$,   \begin{eqnarray}T_{\tau\circ \sigma}=T_\sigma\circ T_\tau,\end{eqnarray} 
 from which follows that \begin{eqnarray}T_\sigma^{-1}=T_{\sigma^{-1}}.\end{eqnarray}
 Now
 for every $\sigma\in H$, 
 any $N$ and $A_1,...,A_{2N+1}\in {\cal F}$, one can see that the following equality holds
 \begin{eqnarray}
\mu(X'_{\sigma(-N)}\in A_1,...,X'_{\sigma(N)}\in A_{2N+1})=\mu(X'_{-N}\in A_1,...,X'_{N}\in A_{2N+1})\hskip 0,8 cm (\textit{exchang}).
\end{eqnarray}
and that it is also equivalent to
\begin{eqnarray}
\mu\circ T_\sigma^{-1}(X'_{-N}\in A_1,...,X'_{N}\in A_{2N+1})=\mu(X'_{-N}\in A_1,...,X'_{N}\in A_{2N+1}).
\end{eqnarray}
In conclusion, the preceding 
   shows that the exchangeability of the process  $(X_n)_{n\ge 1}$ is the same as the exchangeability  of the process $(X'_n)_{n\in \mathbb Z}$ and it is also equivalent to  the invariance of $\mu $, the law of  $\hskip 0,1 cm  (X_n')_{n\in \mathbb Z},\hskip 0,1 cm $ by the transformation $T_\tau$, for all $\tau\in H$, and in particular, implies,  as noted before,  the invariance of $\hskip 0,1 cm \mu\hskip 0,1 cm$ by the shift $ \hskip 0,1 cm S$. Henceforth, the process $(X_n)$ we consider will be indexed by $\mathbb Z$, and furthermore $X_n$ will be  the $n^{\textit{th}}$ coordinate function on $\hskip 0,2 cm \Omega:=K^\mathbb Z$.\\
   
A slight generalization of exchangeability is given by the following\\

{\bf{Definition 2}}\\
\textit{We say that a sequence  $(X_n)_{n\in \mathbb Z}$ of random variables, with law $\mu $,  is quasi exchangeable if  $\mu\circ T_\sigma^{-1}$ is equivalent to $\mu$, for all permutation $\sigma\in H$.}\\

In this case we denote the Radon-Nikodym derivative of $ \hskip 0,1 cm \mu\circ T_\sigma^{-1} \hskip 0,1 cm$ with respect to $ \hskip 0,2 cm \mu, \hskip 0,1 cm$ by  $ \hskip 0,1 cm \phi_\sigma$
\begin{eqnarray}
\phi_\sigma:={d\mu\circ T_\sigma^{-1}\over d\mu}.
\end{eqnarray}
Since we are only interested exclusively in stationary sequences of random variables, the following remark may justify Definition 2.\\
{\bf{Remark 1}}\\
\textit{ Naturally we shall say that a unilateral  sequence $(X_n)_{n\ge 1}$  of random variables is quasi-exchangeable if its law $P$ is equivalent to the law  $P_\sigma$  of the sequence $(X_{\sigma(n)})_{n\ge 1},$ for all permutation $\sigma \in H_1$, where $H_1$ is the group of all finite permutations of the set $\mathbb N$ of natural numbers.\\
Then, if $(X_n)_{n\ge 1}$, with law $P,$ is   quasi-exchangeable  and  stationary, the sequence $(X'_n)_{n\in \mathbb Z}$, with law $\mu$ defined by $(4)$, will be  shift invariant and quasi-exchangeable. }\\

{\bf{Definition 3}}\\
\textit{Let $(X_n)_{n\in \mathbb Z}$ be a quasi exchangeable process, with law $\hskip 0,1 cm \mu.\hskip 0,1 cm$   If the family  $\{{d\mu\circ T_\sigma^{-1}\over d\mu}:\sigma\in H\}$ of all Radon-Nikodym derivatives is uniformly integrable, we say that the process  $X$ is quasi-exchangeable with uniformly integrable densities. }\\

We shall also use the following notations.\\
If  $L$ and $s$ are integers with $L\ge 0$ and $s\ge 1$,  and   $A_{-L},...,A_s$ are  measurable subsets  of $K$, we set 
\begin{eqnarray}
{\Pi}(A_{-L},...,A_0):=\{\omega\in \Omega: \omega_{-L}\in A_{-L},...,\omega_0\in A_0\}\\
F(A_1,...,A_s):=\{\omega \in \Omega:\omega_1\in A_1,...,\omega_s\in A_s\} .
\end{eqnarray}
and for all $I\subset \mathbb Z$,
\begin{eqnarray*}
{\cal A}_I:= \textit{the smallest algebra containing the  sets} \hskip 0,2 cm \{\omega\in \Omega:\omega_j\in A\}, j\in I, A\in {\cal F},\\
\textit{and}\hskip 2 cm 
{\cal B}_I:=\textit{the sigma algebra generated by }\hskip 0,2 cm  {\cal A}_I.
\end{eqnarray*}
In the following particular cases:
\begin{eqnarray*}
I=\{n\in \mathbb Z:n\le 0\}\hskip 0,4 cm  {\cal A}_I \hskip 0,2 cm  \textit{ is denoted}\hskip 0,2 cm  {\cal A}_{\le 0}\\
I=\{n\in \mathbb Z:n\ge  p\}\hskip 0,4 cm  {\cal A}_I \hskip 0,2 cm  \textit{ is denoted}\hskip 0,2 cm  {\cal A}_{\ge p}\\
I=\mathbb Z \hskip 0,4  cm  {\cal A}_I \hskip 0,2 cm  \textit{ is denoted}\hskip 0,2 cm  {\cal A}.
\end{eqnarray*}
The same notation will be used for ${\cal B}_I$, in particular  $\hskip 0,3 cm {\cal B}_{\mathbb Z}={\cal B}$. \\
Similarly, if $\mu$ is a measure on $\Omega$, then $\mu_I$ denotes the restriction of $\mu$ to the sigma-algebra ${\cal B}_I$.  
Also, $\hskip 0,2 cm M_1(Y)\hskip 0,2 cm $ will denote the set of all Radon probability measures on the topological\\ space
$\hskip 0,2 cm Y$, and if $(Z,{\cal G}, m)$ is a probability space and ${\cal G}_1$ is a sub sigma-algebra of ${\cal G}$, 
the conditional expectation of a function $f\in L^1(Z,{\cal G},m)$ given  ${\cal G}_1$ will be denoted by $E_m[f\mid {\cal G}_1]$, or, if there is no confusion on the measure $m$, simply by $E[f\mid{\cal G}_1]$.\\
Also, the smallest sigma-algebra which contains  two sigma-algebras   ${\cal F}_1$ and ${\cal F}_2$  is denoted by $\hskip 0,2 cm {\cal F}_1\vee {\cal F}_2$, and the smallest one containing a family $\{{\cal F}_j:j\in J\}$, of sigma-algebras is denoted $\bigvee_{j\in J} {\cal F}_j$.  The complement of a subset $A$ is denoted $A^c$.
\section{The main results}
Recall that $\hskip 0,1 cm \Omega=K^\mathbb Z,\hskip 0,1 cm$ and $\hskip 0,1 cm {\cal F}\hskip 0,1 cm$ is a sigma algebra of subsets of $\hskip 0,1 cm K.\hskip 0,1 cm$ In the following theorem we suppose that $\hskip 0,1 cm{\cal F}\hskip 0,1 cm$ is separable. \\ 
{\bf{Theorem 1}}\\
\textit{Let $(X_n)_{n\in \mathbb Z}$ be a stationary  quasi exchangeable process, with law $\mu$, such that the dynamical system $(\Omega,S,\mu)$ is ergodic. Then the process is isomorphic to a Bernoulli process.\\
More precisely, if the state space $K$ is finite the system $(\Omega,S,\mu)$ is  "faiblement de Bernoulli"}.\\

{\bf{Proof}}  We consider first the case where  the state space $K$ is finite, and then we prove  that  $(\Omega, S,\mu)$ is "faiblement de Bernoulli". For this, we shall use the ergodicity and the quasi-exchangeability under a particular infinite   family of permutations, in order to find a sequence of measures converging to $\hskip 0,2 cm \mu_{\le 0}\times \mu_{\ge 1}\hskip 0,2 cm$ on all cylinders and such that any measure in the sequence, coincides with $\mu$ on the two-sided tail sigma-algebra $\hskip 0,2 cm {\cal T}:=\bigcap_n ( {\cal B}_{\le -n}\vee {\cal B}_{\ge n})$. The detail follows.\\
For all natural numbers  $ k$ and $ P,\hskip 0,2 cm $ such that $\hskip 0,2 cm  1\le k<P,\hskip 0,2 cm $ let us consider the permutation ( involution) $\hskip 0,2 cm  \sigma:=\sigma_{P,k}\in H\hskip 0,2 cm  $   which translates the "interval" $\hskip 0,1 cm I:=\mathbb N\cap [1,...,k]\hskip 0,1 cm $  by $\hskip 0,1 cm P$,  translates the $^"$interval$^"$ $\hskip 0,1 cm P+I=\mathbb N\cap [P+1,...,P+k]\hskip 0,1 cm $ by $\hskip 0,1 cm -P,\hskip 0,1 cm$ and leaves fixed all $\hskip 0,1 cm n\in \mathbb Z, \hskip 0,1 cm$ which are not in the disjoint union  $\hskip 0,1 cm  I\cup (P+I)\hskip 0,1 cm$, that is, which is   defined by
\begin{eqnarray*}
\sigma(j)=P+j, \hskip 0,2 cm \textit{and}  \hskip 0,2 cm  
 \sigma(P+j)=j,\hskip 0,2 cm \textit{if} \hskip 0,2 cm 1\le j\le k\\
\sigma(n)=n\hskip 0,2 cm \textit{if} \hskip 0,2 cm n\notin\{1,...,k\}\cup \{P+j:j=1,...,k\},
\end{eqnarray*}
so that
\begin{eqnarray*}
 T_\sigma(\omega)_n=\omega_n , \hskip 0,1 cm \textit{if}\hskip 0,2 cm n\notin\{1,...,k\}\cup \{P+j:j=1,...,k\}\\
T_\sigma(\omega)_j=\omega_{P+j}\hskip 0,5 cm  \textit{and} \hskip 0,5 cm  T_\sigma(\omega)_{P+j}=\omega_j  \hskip 0,5 cm \textit{for }\hskip 0,5 cm 1\le j\le k .
\end{eqnarray*}

Recall that  $\hskip 0,1 cm  \phi_\sigma:={d\mu\circ T_\sigma^{-1}\over d\mu},\hskip 0,1 cm $ and for all
   natural numbers  $\hskip 0,2 cm  L\ge 0 $ and $n \ge 1$,  consider    any measurable subsets $\hskip 0,2 cm  A_{-L},...,A_k, \hskip 0,2 cm $  $B_1,...,B_k$ and $\hskip 0,2 cm  C_1,...,C_n\hskip 0,2 cm $  of $\hskip 0,1 cm K,\hskip 0,1 cm$ and any $\hskip 0,3 cm T\in{\cal T}.\hskip 0,3 cm$ Using the notation  as in (10) and (11), and setting
   \begin{eqnarray}
E=\Pi(A_{-L},...,A_0):=\{\omega\in \Omega: \omega_{-L}\in A_{-L},...,\omega_0\in A_0\} ,
\end{eqnarray}

 then, for $\hskip 0,1 cm n+k<P,\hskip 0,1 cm$ the equality
\begin{eqnarray*}
\mu(T_\sigma^{-1}[E\cap F(A_1,...,A_k)\cap S^{-k} F(C_1,...,C_n)\cap S^{-P}F(B_1,...,B_k)\cap T])=\\\mu(E\cap F(B_1,...,B_k)\cap S^{-k}F(C_1,...,C_n)\cap S^{-P}F(A_1,...,A_k)\cap T)
\end{eqnarray*}
which holds by the definitions of $ \sigma$ and $T_\sigma$, 
reads also, in view of quasi-invariance,  as
\begin{eqnarray}
\int _{E\cap F(A_1,...,A_k)\cap S^{-k} F(C_1,...,C_n)\cap T}1_{F(B_1,...,B_k)}\circ S^P\phi_\sigma d\mu= \nonumber \\
\int_{E\cap F(B_1,...,B_k)\cap S^{-k}F(C_1,...,C_n)\cap T}1_{ F(A_1,...,A_k)}\circ S^Pd\mu .
\end{eqnarray}
For every integer $\hskip 0,1 cm Q\ge R+k+1,\hskip 0,1 cm$ where $ R$ is a fixed integer, $\hskip 0,1 cm R\ge n,\hskip 0,1 cm$  consider  the following two functions $\xi_{B_1,...,B_k}^Q\in L^1(\mu)$ and $\psi_{A_1,...,A_k}^Q\in L^\infty(\mu)$, defined by 
\begin{eqnarray*}
\xi_{B_1,...,B_k}^Q:={1\over Q}\sum_{P=R+k+1}^Q 1_{F(B_1,...,B_k)}\circ S^P\phi_{\sigma_{P,k}},\\
\textit{and}\hskip 0,5 cm \psi_{A_1,...,A_k}^Q:={1\over Q}\sum_{P=R+k+1}^Q 1_{F(A_1,...,A_k)}\circ S^P.
\end{eqnarray*}
Then, by $(13)$, we obtain 
\begin{eqnarray*}
\int _{E\cap F(A_1,...,A_k)\cap S^{-k} F(C_1,...,C_n)\cap T}\xi_{B_1,...,B_k}^Q d\mu=\\
\nonumber \int _{E\cap F(B_1,...,B_k)\cap S^{-k} F(C_1,...,C_n)\cap T}\psi_{A_1,...,A_k}^Q d\mu, 
\end{eqnarray*}
which we write  as
\begin{eqnarray*}
\int _{E\cap S^{-k} F(C_1,...,C_n)\cap T} 1_{F(A_1,...,A_k)}\xi_{B_1,...,B_k}^Q d\mu=\\
\nonumber \int _{E\cap S^{-k} F(C_1,...,C_n)\cap T}1_{F(B_1,...,B_k)}\psi_{A_1,...,A_k}^Q d\mu.
\end{eqnarray*}
This last equality,  holding  true for all $L\ge 1,\hskip 0,1 cm $ $1\le n\le R$, all $E$ as in (12),  all $C_1,...,C_n$, and all $T\in{\cal T}$,   means that
\begin{eqnarray}
E[ 1_{F(A_1,...,A_k)}\xi_{B_1,...,B_k}^Q\mid {\cal M}_R]=E[ 1_{F(B_1,...,B_k)}\psi_{A_1,...,A_k}^Q\mid {\cal M}_R], 
\end{eqnarray}
where ${\cal M}_R$ is the sigma-algebra
\begin{eqnarray*}
{\cal M}_R:={\cal B}_{\le 0}\vee {\cal B}_{\{k+1,...,k+1+R\}}\vee {\cal T}.
\end{eqnarray*}

Now, $\hskip 0,1 cm k\hskip 0,1 cm$ being fixed,  by ergodicity,  the sequence $\hskip 0,1 cm (\psi_{A_1,...,A_k}^Q)_{Q\ge R+k+1}\hskip 0,1 cm$ converges in $L^2(\mu)$ norm to the constant $\mu(F(A_1,...,A_k))$. Then, by Cauchy Schwartz for example,  the sequence $\hskip 0,1 cm (1_{F(B_1,...,B_k)} \psi_{A_1,...,A_k}^Q)_{Q\ge R+k+1}\hskip 0,1 cm $ converges in $\hskip 0,1 cm L^2(\mu)\hskip 0,1 cm $ norm also to $\mu(F(A_1,...,A_k)) 1_{F(B_1,...,B_k)}$. It follows that the sequence in the left side of $(14)$  converges in norm $\hskip 0,1 cm L^2(\mu),\hskip 0,1 cm$ \\ to  $\hskip 0,1 cm \mu(F(A_1,...,A_k)) E[1_{F(B_1,...,B_k)}\mid {\cal M}_R]$. \\

In particular,
\begin{eqnarray}
\lim_Q \sup_{M\in {\cal M}_R } \mid \int_{M\cap F(A_1,...,A_k)} \xi_{B_1...,B_k}^Qd\mu -\mu(M\cap F(B_1,...,B_k)) \mu(F(A_1,...,A_k))\mid =0. 
\end{eqnarray}
The sequence $\hskip 0,2 cm(\xi_{(B_1,...,B_k)}^Q)_{Q\ge R+k+1},\hskip 0,2 cm$ is bounded in  $\hskip 0,2 cm L^1(\mu):\hskip 0,2 cm$\begin{eqnarray*}
\mid\mid \xi_{(B_1,...B_k)}^Q\mid\mid_1\le \mid \mid {1\over Q}\sum_{P=R+k+1}^Q \phi_{\sigma_{P,k}}\mid\mid_1={Q-R\over Q}\le 1.
\end{eqnarray*}\\
 Then it is bounded in the bidual $\hskip 0,2 cm L^{\infty *}(\mu)\hskip 0,2 cm $ of $\hskip 0,2 cm L^1(\mu).\hskip 0,2 cm $  Hence, by Alaoglu-Bourbaki Theorem( [6], Theorem 2, p.424), this sequence has  at least one weak-star cluster point. Let $\hskip 0,1 cm \eta_{(B_1,...,B_k)}^k\in L^{\infty*}(\mu)\hskip 0,1 cm $ be such a cluster point. Then $\hskip 0,2 cm \eta_{B_1,...,B_k}^k\hskip 0,2 cm$ is positive,  because for each $Q$, $\hskip 0,2 cm \xi_{(B_1,...,B_k)}^Q\hskip 0,2 cm $  is positive.  Also, for any $x^*\in L^\infty(\mu)$,  there exists a subsequence of natural numbers  $(Q_j)_{j\ge 1}=(Q_j(x^*))_{j\ge 1}$, which may depend on $x^*$, converging to infinity such that
\begin{eqnarray*}
\eta^k_{(B_1,...,B_k)}(x^*)=\lim_j \xi_{B_1,...,B_k}^{Q_j(x^*)}(x^*).
\end{eqnarray*}
Taking, in particular,  $\hskip 0,1 cm x^*=x_0^*,\hskip 0,1 cm$ where $\hskip 0,1 cm x_0^*=1_{M\cap F(A_1,...,A_k)}$, we obtain
\begin{eqnarray}
\eta_{(B_1,...,B_k)}^k(M\cap F(A_1,...,A_k))=\lim_j \int_{M\cap F(A_1,...,A_k)} \xi_{B_1,...,B_k}^{Q_j(x_0^*)}d\mu. 
\end{eqnarray}
Then, in view of $(15)$, the limit in $(16)$ is independent of the sequence $(Q_j(x_0^*))$, and the following equality holds
\begin{eqnarray}
\eta_{B_1,...,B_k} ^k(M\cap F(A_1,...,A_k))=\mu(M\cap F(B_1,...,B_k))\mu(F(A_1,...,A_k))
\end{eqnarray}
for all $M\in {\cal M}_R$.\\
In the particular case, where  $\hskip 0,2 cm B_1=...=B_k=K,\hskip 0,2 cm$ let us denote $\hskip 0,2 cm \eta_{B_1,...,B_k}^k\hskip 0,2 cm$ simply by $\hskip 0,2 cm \eta^k.\hskip 0,2 cm$ Then, by $(17)$,
\begin{eqnarray}
\eta^k(M\cap F(A_1,...,A_k))=\mu(M)\mu(F(A_1,...,A_k)), \forall  M\in {\cal M}_R,
\end{eqnarray}
which means that, under  $\eta^k$, the algebra ${\cal A}_{\{1,...,k\}}$ and the sigma algebra ${\cal M}_R$ are independent( Note that as element of $L^{\infty *}(\mu)$, $\eta^k$ may have a non null purely finitely additive part).
  Recall that 
 \begin{eqnarray*}
 {\cal T}:=\bigcap_{n\ge 1}({\cal B}_{\le -n}\vee {\cal B}_{\ge n}),
\end{eqnarray*}
and  set \begin{eqnarray*}\hskip 0,1 cm {\cal M}={\cal T}\vee {\cal B}_{\le 0}.\hskip 0,1 cm\end{eqnarray*} Then $\hskip 0,1 cm {\cal M}\subset {\cal M}_R, \hskip 0,1 cm$
 and, by $(18)$,  we have in particular
\begin{eqnarray}
\eta^k(A\cap M)= \mu(A)\mu(M), \hskip 0,1 cm \forall A\in{\cal A}_{\{1,...,k\}},\hskip 0,1 cm \forall M\in {\cal M}.
\end{eqnarray}
Also, if we  denote by $\hskip 0,1 cm {\cal D}_k,\hskip 0,1 cm$ the smallest algebra containing  $\hskip 0,1 cm{\cal M}\cup  {\cal A}_{\{1,...,k\}},\hskip 0,1 cm$ then $\hskip 0,1 cm{\cal D}_k\subset {\cal D}_{k+1},\hskip 0,1 cm$ 
and $(19)$ implies that $ \hskip 0,1 cm \eta^{k+1}\hskip 0,1 cm$ extends $\hskip 0,1 cm \eta^k\hskip 0,1 cm$ from $\hskip 0,1 cm {\cal D}_k\hskip 0,1 cm$ to $\hskip 0,1 cm {\cal D}_{k+1}$. By induction we then have, for every $n\ge 0$, 
\begin{eqnarray}
\eta^{k+n}(A\cap M)= \mu(A)\mu(M), \hskip 0,1 cm \forall A\in{\cal A}_{\{1,...,k\}},\hskip 0,1 cm \forall M\in {\cal M}.
\end{eqnarray}

Now, the sequence $\hskip 0,2 cm (\eta^k)_{k\ge 1}\hskip 0,2 cm $ is norm bounded in $\hskip 0,2 cm L^{\infty*}(\mu),\hskip 0,2 cm$ since by positivity and by $(18)$ we have
\begin{eqnarray*}
\mid\mid \eta^k \mid\mid_{L^{\infty *}(\mu)}=\eta^k(\Omega)=1.
\end{eqnarray*}
Then $\hskip 0,1 cm (\eta^k)_{k\ge 1}\hskip 0,1 cm$  has at least one weak-star cluster point. By $(20)$, every such cluster point, say $\hskip 0,2 cm \eta,\hskip 0,2 cm$ verifies
\begin{eqnarray}
\eta(A\cap M)= \mu(A)\mu(M), \hskip 0,1 cm \forall A\in{\cal A}_{\{1,...,k\}},\hskip 0,1 cm \forall k\ge 1,\hskip 0,1 cm \forall M\in {\cal M}.
\end{eqnarray}
which means
\begin{eqnarray*}
\eta(A\cap M)= \mu(A)\mu(M), \hskip 0,1 cm \forall A\in{\cal A}_{\ge 1},\hskip 0,1 cm \hskip 0,1 cm \forall M\in {\cal M}.
\end{eqnarray*}
In particular
\begin{eqnarray}
\eta(A\cap B)=\mu(A)\mu(B),\hskip 0,2 cm \forall A\in {\cal B}_{\le 0},\hskip 0,2 cm \forall B\in {\cal A}_{\ge 1}.
\end{eqnarray}
Notice that $(22)$   implies that $  \hskip 0,2 cm \eta \hskip 0,2 cm $ and $  \hskip 0,2 cm \mu_{\le 0}\times \mu_{\ge 1}  \hskip 0,2 cm$ coincide on all cylinders. We recall that, for any subset $\hskip 0,2 cm I\subset \mathbb Z$,  $\hskip 0,2 cm \mu_I$ denotes the restriction of $\mu$ to the sigma algebra ${\cal B}_I$. 
Also $(21)$ implies that $ \hskip 0,2 cm \eta \hskip 0,2 cm$ and $ \hskip 0,2 cm \mu \hskip 0,2 cm$ coincide on $ \hskip 0,2 cm {\cal T}.$ \\
It follows that if $\hskip 0,1 cm {\cal L}\hskip 0,1 cm$ is the algebra generated by $\hskip 0,1 cm  {\cal M}\cup {\cal A}_{\ge 1}$, then clearly $\hskip 0,1 cm {\cal L}\hskip 0,1 cm$ contains $\hskip 0,1 cm {\cal A}\hskip 0,1 cm$ and, by $(21)$,   that $\hskip 0,1 cm \eta\hskip 0,1 cm$ and  $\hskip 0,1 cm \mu_{_{\le 0}} \times \mu_{_{\ge 1}}\hskip 0,1 cm$ coincide also on $\hskip 0,1 cm {\cal L},\hskip 0,1 cm$  
so  that $\hskip 0,1 cm \eta\hskip 0,1 cm$ is countably additive on  $\hskip 0,1 cm {\cal L}.\hskip 0,1 cm$  We also note  that $\hskip 0,1 cm {\cal B}\hskip 0,1 cm$ is the smallest sigma-algebra containing $\hskip 0,1 cm {\cal L}.\hskip 0,1 cm$ Then the unique countably additive measure $\hskip 0,1 cm \tilde{\eta}\hskip 0,1 cm$ extending $\hskip 0,1 cm \eta\hskip 0,1 cm$ to the sigma algebra $\hskip 0,1 cm {\cal B},\hskip 0,1 cm$  is the measure $\hskip 0,1 cm \mu_{\le 0}\times \mu_{\ge 1}.\hskip 0,1 cm$\\
Then, in particular,  $\hskip 0,2 cm \tilde{\eta}\hskip 0,2 cm$   verifies 
\begin{eqnarray}
\tilde{\eta}(A\cap B):=\mu(A)\mu(B), \hskip 0,2 cm A\in{\cal B}_{\le 0}, \hskip 0,2 cm B\in {\cal B}_{\ge 1}. 
\end{eqnarray}
Since $ {\cal T}\subset  {\cal M}$, we also have, by $(21)$, as noted before,
\begin{eqnarray}
\tilde{\eta}(M)=\eta(M)=\mu(M),\hskip 0,2 cm \forall M\in {\cal T}.
\end{eqnarray}
In conclusion, the countably additive measure $\hskip 0,1 cm \tilde{\eta}\hskip 0,1 cm$ on $\hskip 0,1 cm {\cal B}\hskip 0,1 cm$ satisfies the equalities $(23)$ and $(24)$, which means, when the state space $K$ is finite, that the system $(\Omega,S,\mu)$ is "faiblement de Bernoulli" ( [11], D\' efinition 3). Since in this case, "faiblement de  Bernoulli" is equivalent to weak Bernoulli ( [11], Proposition 2) and also, a system which is weak Bernoulli   is isomorphic to a Bernoulli system [7], the proof is complete in this finite case.\\
To end the proof in the general case, we establish first the following Lemma, and we appeal to a result of [17] (Theorem 2), saying that if   $T$ is a one to one invertible measure preserving transformation on a measure algebra ${\cal L }$, if ${\cal L}$ is the increasing union of invariant sub-sigma algebras ${\cal L}_i$ such that $T$ restricted to each ${\cal L}_i$ is a Bernoulli shift, then $T$ is a generalized Bernoulli shift.\\

{\bf{Lemma 1}}\\
\textit{ Let $\Omega=K^\mathbb Z$, $S$ the shift on $\Omega$ and $\mu$ be an $S$-invariant probability measure on $\Omega$. Let $X_n$ be the $n^\textit{th}$ coordinate function on $\Omega$ 
and  ${\cal P}=\{A_0,...,A_{k-1}\}$ be a finite measurable partition of $K$.\\ Let $\Omega_1:=\{0,...,k-1\}^\mathbb Z,$ $S_1$ the shift on $\Omega_1$, $Y_n$ the $n^\textit{th}$ coordinate function on $\Omega_1$  and $\theta:\Omega\rightarrow \Omega_1$ the factor map defined by:
\begin{eqnarray*}
\theta(x)=y\iff (S^nx)_0\in A_{y(n)}, \forall n\in \mathbb Z\iff x_n\in A_{y(n)},\forall n\in \mathbb Z,
\end{eqnarray*}
for any $x=(x_n)_{n\in \mathbb Z}\in \Omega$, 
that is
\begin{eqnarray*}
\theta(x)(n)=j\iff (S^nx)_0\in A_j\iff x_n\in A_j.
\end{eqnarray*}
Let, for all $n$, 
 \begin{eqnarray*}
Z_n:=Y_n\circ \theta.
\end{eqnarray*}
Then \\
(1) the quasi-exchangeability[ respectively exchangeability] of the process $X=(X_n)$ implies the same property for the process $Z=(Z_n)$.\\
(2) the ergodicity of the process $(X_n)$ implies the ergodicity of $(Z_n)$.}\\ 

{\bf{Proof}} Let $\nu:=\mu\circ \theta^{-1}$ and, for any permutation $\tau$,  $R_\tau:\Omega_1\rightarrow \Omega_1,$ the transformation defined by
\begin{eqnarray*}
R_\tau(y)(n)=y(\tau(n)), \forall n\in \mathbb Z,\forall y\in \Omega_1.
\end{eqnarray*}
Then  the quasi-exchangeability[ respectively exchangeability] of $(Z_n)$  follows from the commutation relationship
\begin{eqnarray*}
R_\sigma\circ \theta=\theta\circ T_\sigma
\end{eqnarray*}
because if this relationship  holds we will get
\begin{eqnarray*}
\nu\circ R_\sigma^{-1}=\mu\circ \theta^{-1}\circ R_\sigma^{-1}=\mu\circ T_\sigma^{-1}\circ \theta^{-1}
\end{eqnarray*}
and then
\begin{eqnarray*}
\nu\circ R_\sigma^{-1}(A)=0\iff (\mu\circ T_\sigma^{-1})(\theta^{-1}A)=0\iff \mu(\theta^{-1}A)=0\iff \nu(A)=0.
\end{eqnarray*}
[Respectively
\begin{eqnarray*}
\nu\circ R_\sigma^{-1}(A)= (\mu\circ T_\sigma^{-1})(\theta^{-1}A)=\mu(\theta^{-1}A)=\nu(A).]
\end{eqnarray*}

The commutation relation is a consequence of    the following
\begin{eqnarray*}
R_\sigma(\theta(x))(n)=j\iff \theta(x)(\sigma(n))=j\iff x_{\sigma(n)}\in A_j\\
\iff T_\sigma(x)(n)\in A_j\iff \theta(T_\sigma(x))(n)=j
\end{eqnarray*}
where $x\in \Omega, n\in \mathbb Z$, and $j=0,1,...,k-1$.\\

{\bf{Corollary 1}}\\
\textit{Under the same conditions of the above Lemma, if the process $X$ with law $\mu$ is quasi-exchangeable and the system $(\Omega,S,\mu)$ is ergodic, then the system $(\Omega_1,S_1,\mu\circ \theta^{-1})$ is isomorphic to a Bernoulli system.}\\

From this Corollary, it follows that if ${\cal P} $ is a finite measurable partition of $K$ and  if $\hskip 0,1 cm {\cal P}_S:=\bigvee_{n\in \mathbb Z} S^n{\cal P},\hskip 0,1 cm$ is the smallest sigma-algebra containing $\hskip 0,1 cm \cup_{n\in \mathbb Z} S^n{\cal P},\hskip 0,1 cm$  then the system $(\Omega, S,{\cal P}_S)$ is isomorphic to a Bernoulli system and is a factor of $(\Omega,S,{\cal B}(\Omega))$.\\Now, if ${\cal Q}$ is a finite partition of $K$, which refines ${\cal P}$, then ${\cal P}_S\subset {\cal Q}_S$, that is   $(\Omega, {\cal P}_S, \mu)$ is a factor of    $(\Omega, {\cal Q}_S, \mu)$.  Then,  because  ${\cal B}$ is generated by $\cup\{ {\cal P}_S: {\cal P} \textit{ is a finite partition of}\hskip 0,1 cm  K\}$,  and since this last union can also be written as an  increasing union $\hskip 0,1 cm \cup_i {\cal P}_S^i ,\hskip 0,1 cm$ where for each $\hskip 0,1 cm i,\hskip 0,1 cm$ $(\Omega, S,{\cal P}^i_S)\hskip 0,1 cm$ is isomorphic to a Bernoulli shift,   we conclude, by the above mentioned Theorem of [17], that  the system $(\Omega,S,{\cal B},\mu)$ is a generalized Bernoulli system,  and the proof of Theorem 1 is achieved.$\square$\\

Recall  that the exchangeability means the equalities $\mu\circ T_\sigma^{-1}=\mu$, for all $\sigma\in H$,  so that the following result generalizes De Finetti's Theorem. \\

{\bf{Theorem 2}}\\
\textit{Let $X=(X_n)_{n\in \mathbb Z}$ be a stationary process. Then the following properties are equivalent\\\
(1)  $X$ is exchangeable.\\
(2) $X$ is quasi-exchangeable with uniformly integrable densities.\\
(3) $X$ is a mixture of Bernoulli processes.}

{\bf{Proof}} The equivalence $(1)\iff (3)$ is the generalization by Hewitt-Savage of De Finetti's Theorem to presentable spaces. The implication  $(1)\Rightarrow (2)$ is trivial. We show that $(2)$ implies $(3)$.  For this, let $\mu$ denotes the law of the process, and  let $H_N:=\{\sigma\in H: \sigma(n)=n, \forall n, \mid n\mid >N\}$. Then $H_N$ is a subgroup of $H$ and $H=\cup_N H_N$. For every $N\ge 1$, consider  the measure 
\begin{eqnarray*}
\nu_N:={1\over \textit{card}(H_N)}\sum_{\sigma\in H_N} \mu\circ T_{\sigma}^{-1}=({1\over \textit{card}(H_N)}\sum_{\sigma\in H_N}  \phi_\sigma) \mu,
\end{eqnarray*}
where $\textit{card}(H_N)=(2N+1)!$ denotes the cardinality of $H_N$.\\
Scheme of the proof: We prove first that any cluster point $\nu$ of the sequence $(\nu_N)_{N\ge 1}$ is invariant by $T_\sigma$, for any permutation $\sigma\in H$ (and thus the Hewitt-Savage generalization of De Finetti's Theorem applies to that cluster point $\nu$). The hypothesis implies that $\nu$ is absolutely continuous with respect to $ \mu$, and using the definition of $\nu$, we prove that also $\mu$ is absolutely continuous with respect to $\nu$ and we conclude using the Hewitt-Savage mentioned Theorem( or also by Corollary 2 and Remark 3 below).\\   
The details are as  follows.\\
Notice that the uniform integrability of the family $\{\phi_\sigma:\sigma\in H\}$ implies the same property for the family
$\{{1\over \textit{card}(H_N)}\sum_{\sigma\in H_N}  \phi_\sigma:N\ge 1\}$, so that,  if $\hskip 0,2 cm f_N:={1\over \textit{card}(H_N)}\sum_{\sigma\in H_N}  \phi_\sigma$,  then, according to Dunford-Pettis Theorem,  $\hskip 0,2 cm (f_N)$  admits a subsequence $\hskip 0,2 cm (f_{N_k}),\hskip 0,2 cm$ which converges weakly in $\hskip 0,2 cm L^1(\mu),\hskip 0,2 cm$ to some $f\in L^1(\mu).\hskip 0,2 cm$
It follows that the sequence of probability measures $\nu_{N_k}$  converges to $\nu:=f\mu$, in the sense that
\begin{eqnarray*}
\forall A\in {\cal B}, \hskip 0,2 cm \nu(A)=\lim_k \nu_{N_k}(A).
\end{eqnarray*}
 Now, because, for every $ \hskip 0,1 cm \tau\in H,\hskip 0,1 cm$ there is $k_0(\tau)$ such that $\tau\in H_{N_k} $, for any $k\ge k_0(\tau)$,   the measure $\nu$ is invariant by $T_\tau$, for all $\tau\in H$. In fact let $A$ be a measurable set, and $\tau\in H$. Then, with $M:=N_k$, and using (5), the following equalities hold
\begin{eqnarray*}
\nu_{N_k}(T_\tau^{-1}A)={1\over \textit{card}(H_M)}\sum_{\sigma\in H_M} \mu\circ T_{\sigma}^{-1}(T_\tau^{-1}A)={1\over \textit{card}(H_M)}\sum_{\sigma\in H_M} \mu\circ (T_\tau T_{\sigma})^{-1}(A)\\
={1\over \textit{card}(H_M)}\sum_{\sigma\in H_M} \mu\circ T_{\sigma\circ \tau}^{-1}(A)={1\over \textit{card}(H_M)}\sum_{\sigma\in (H_M)\tau} \mu\circ T_{\sigma}^{-1}(A).
\end{eqnarray*}
But, as noted before,  for $k$ big enough, $\hskip 0,1 cm \tau\in H_{N_k}=H_M,\hskip 0,1 cm$ so that $\hskip 0,1 cm H_M=(H_M)\tau,\hskip 0,1 cm$ since $\hskip 0,1 cm H_M\hskip 0,1 cm$ is a group, and then
\begin{eqnarray*}
\nu_{N_k}(T_\tau^{-1}A)={1\over \textit{card}(H_M)}\sum_{\sigma\in (H_M)\tau} \mu\circ T_{\sigma}^{-1}(A)={1\over \textit{card}(H_M)}\sum_{\sigma\in H_M} \mu\circ T_{\sigma}^{-1}(A) =\nu_{N_k}(A)
\end{eqnarray*}
and consequently  $\hskip 0,1 cm \nu(T_\tau^{-1}(A))=\nu(A)$.

It follows that $\nu$ is also invariant by the shift $S$, (hence it will be  equal $\mu$, if  $\mu$ was ergodic, and  hus $\mu$ will be  invariant by $H$, hence $\mu$ will be  Bernoulli, by the Hewitt-Savage generalization of De Finetti's Theorem, or by Proposition 1 below).\\
\\Suppose now that the system is not necessarily ergodic. By a slight adaptation of the Hewitt-Savage generalization of De Finetti's Theorem [8] or also, in the case where $K$ is compact,  by Corollary 2 and Remark 3 below,  being invariant by all $T_\tau$, the measure $\nu$ is given by an average of independent measures:
\begin{eqnarray}
\nonumber \textit{there exists a probability measure}  \hskip 0,5 cm \beta \in M_1(\tilde{K}),\hskip 0,5 cm \textit{such that}\\
\nu(A)=\int_{\tilde{K}} (\tilde{\pi}(A)) d\beta(\tilde{\pi})=\int_{\tilde{K}}( \int_ {K^\mathbb Z} 1_A d\tilde{\pi}) d\beta(\tilde{\pi})
\end{eqnarray}
for all $A\in {\cal B}_a(K^{\mathbb Z})$, the Baire sigma algebra of $\hskip 0,2 cm K^\mathbb Z,\hskip 0,2 cm $ where $\tilde{K}=M_1(K^\mathbb Z)$ and where,  for any probability measure $\pi$ on $K$, the probability measure $\hskip 0,1 cm  \tilde{\pi}\hskip 0,1 cm $  denotes   the corresponding  product measure on $\hskip 0,2 cm  K^\mathbb Z$: $\hskip 0,2 cm   \tilde{\pi}:=\pi^{\otimes \mathbb Z}$.\\  

By repeated use of Lebesgue dominated convergence theorem and monotone convergence theorem, we deduce from $(25)$  that, 
for any $ {\cal B}_a(K^{\mathbb Z})$ measurable and $\nu$-integrable  $h$,
\begin{eqnarray}
\nu(h)=\int_{\tilde{K}} d\beta(\tilde{\pi})\int _{K^\mathbb Z} h d\tilde{\pi}.
\end{eqnarray}

 Also, $\hskip 0,1 cm \nu=f\mu\hskip 0,1 cm $ is absolutely continuous with repect to $\hskip 0,1 cm \mu,\hskip 0,1 cm$ and thus, $f$ is  $S$-invariant $\mu$ almost everywhere.\\  Now the following implications, which hold for all $\sigma\in H$, 
\begin{eqnarray*}
\nu\circ T_{\sigma}^{-1}=\nu \iff f\circ T_{\sigma}^{-1}\mu\circ T_{\sigma}^{-1}=f\mu\iff f\circ T_\sigma^{-1}\phi_\sigma \mu =f\mu\\ 
\iff f\circ T_\sigma^{-1}\phi_\sigma=f \hskip 0,4 cm \mu \hskip 0,1 cm a.e.\iff f\phi_\sigma \circ T_\sigma=f\circ T_\sigma \hskip 0,4 cm \mu \hskip 0,1 cm \hskip 0,1 cm a.e.\end{eqnarray*}
(the last implications  use the equivalence of $\hskip 0,2 cm  \mu,\hskip 0,2 cm  $ $\hskip 0,2 cm  \mu\circ T_\sigma^{-1}\hskip 0,2 cm  $ and $\hskip 0,2 cm  \mu\circ T_\sigma$)
 implie that
 the set $\hskip 0,1 cm A_0:=\{f=0\}\hskip 0,1 cm $ is mod($\mu$) invariant by $\hskip 0,1 cm  T_\sigma$, for all $\hskip 0,2 cm  \sigma$.
 Then, in particular $$\hskip 0,1 cm \mu(T_\sigma^{-1}(A_0))=\mu(A_0),\hskip 0,1 cm$$ so that
 \begin{eqnarray*}
0=\nu(A_0)=\lim_k \nu_{N_k}(A_0)=\lim_k {1\over \textit{card}(H_{N_k})} \sum_{\sigma\in H_{N_k} }  \mu(T_\sigma^{-1}(A_0))=\mu(A_0).
\end{eqnarray*}
Then $\hskip 0,1 cm \nu\hskip 0,1 cm$ is equivalent to $\hskip 0,1 cm \mu$ so that, for some $S$-invariant  $g\in L^1(\nu)$, $ \mu=g\nu$.\\
Now, according to Lemma 2  below, for some  $E$  with $\nu(E)=1$,  $\hskip 0,1 cm g=g_1$ on $E$ and $g_1$ is  $ {\cal B}_a(K^\mathbb Z)$ measurable.  Then, setting  $\hskip 0,1 cm F:=\cap_{n\in \mathbb Z} S^n E,\hskip 0,1 cm$ we obtain    $\hskip 0,1 cm g=g\circ S^n=g_1\circ S^n=g_1\hskip 0,1 cm $ on $\hskip 0,1 cm F.\hskip 0,1 cm$ We note that $\hskip 0,1 cm g_1\circ S^n\hskip 0,1 cm$ is $\hskip 0,1 cm {\cal B}_a(K^\mathbb Z)\hskip 0,1 cm $ measurable for every $\hskip 0,1 cm n.\hskip 0,1 cm$  Then, for any  $\hskip 0,1 cm A\in {\cal B}_a(K^\mathbb Z), \hskip 0,1 cm$  and any $n
$,
\begin{eqnarray*}
\nu(g1_A)=\nu(g_11_A1_F)=\nu(g_1\circ S^n1_A1_F)=\nu(g_1\circ S^n 1_A)\hskip 0,2 cm \textit{since}\hskip 0,2 cm g_11_F=g_1\circ S^n1_F\\
\textit{but by }(26), \hskip 0,4 cm \nu(g_1\circ S^n 1_A)=\int_{\tilde {K}} d\beta(\tilde{\pi}) \int_{K^\mathbb Z} g_1\circ S^n 1_Ad \tilde{\pi},
\end{eqnarray*}
and then, since this last equality holds for every $n$, we have
\begin{eqnarray*}
\nonumber \nu(g1_A)=\lim_n \int_{\tilde {K}} d\beta(\tilde{\pi}) \int_{K^\mathbb Z} g_1\circ S^n 1_Ad \tilde{\pi}.
\end{eqnarray*}
Then, using  Lebesgue and the mixing property  of the Bernoulli system $\hskip 0,1 cm (K^\mathbb Z,S,\tilde{\pi}),\hskip 0,1 cm$ we obtain
\begin{eqnarray}
\nu(g1_A)=\int_{\tilde{K}} d\beta(\tilde{\pi})\int_{K^\mathbb Z} g_1d \tilde{\pi}\int_{K^\mathbb Z} 1_Ad \tilde{\pi}, 
\end{eqnarray}
because
\begin{eqnarray*}
v_n(\tilde{\pi}):=\int_{K^\mathbb Z} g_1\circ S^n1_Ad\tilde{\pi}\rightarrow \int_{K^\mathbb Z} g_1d\tilde{\pi}\int_{K^\mathbb Z} 1_Ad\tilde{\pi}
\end{eqnarray*}
 and 
 \begin{eqnarray*}
\mid v_n(\tilde{\pi})\mid \le \int_{K^\mathbb Z} g_1d\tilde{\pi},\hskip 2 cm \int_{\tilde{K} }  d\beta(\tilde{\pi})\int_{K^\mathbb Z}  g_1d\tilde{\pi}=\nu(g_1)=1.
\end{eqnarray*}

Now setting $ \hskip 0,1 cm v(\tilde{\pi}):=\int_{K^\mathbb Z}  g_1d\tilde{\pi},\hskip 0,1 cm$ and  $\hskip 0,1 cm \beta_1=v\beta,\hskip 0,1 cm$  the equality $ (27) $ reads
\begin{eqnarray*}
\nu(g1_A)=\int_{\tilde {K}} d\beta(\tilde{\pi}) v(\tilde{\pi})\int_{K^\mathbb Z} 1_Ad \tilde{\pi}, 
\end{eqnarray*}
that is    
\begin{eqnarray*}
\mu(A)=\nu(g1_A)=\int_{\tilde {K}} d\beta_1(\tilde{\pi}) \int_{K^\mathbb Z} 1_Ad \tilde{\pi}=\int_{\tilde {K}}  \tilde{\pi}(A) d\beta_1(\tilde{\pi})
\end{eqnarray*}
and this ends the proof, because
\begin{eqnarray*}
\int_{\tilde {K}} d\beta_1(\tilde{\pi}) =\int_{\tilde {K}} d\beta_1(\tilde{\pi}) \tilde{\pi}(1)=\nu(g)=\mu(1)=1.\square
\end{eqnarray*}

  The following Lemma 2  was used in the proof of Theorem 2. We leave its proof to the reader.\\

{\bf{Lemma 2}}( [9], ex. (12,63), p. 186)\\
\textit{Let $(X,{\cal B}, \nu)$ be a measure space and $f$ be a $\nu$-measurable function. Then\\
there exists a ${\cal B}$-measurable function $g$, which is equal to $f$, $\nu$ almost everywhere.  } \\
In particular\\
\textit{ If $\mu$ is a Borel (respectively  Baire) measure on a topological space and if $f\in L^1(\mu)$ then there exists a Borel( respectively Baire) function $g$ such that $f=g\hskip 0,2 cm $  $\mu$ almost everywhere.}\\

The simplicity of the proof together with  the absence of topological hypothesis on the state space $K$, justifie the following particular case of De Finetti's theorem. \\

{\bf{Proposition 1}}\\
\textit{If $\mu$ is invariant by $T_{\sigma}$ for all $\sigma\in H$, and ergodic for the shift $S$,  then $\mu$ is the product measure.}\\

{\bf{Proof}}
 Let $N\ge 1$, and  for  any $\hskip 0,1 cm P> N+1,\hskip 0,1 cm$  let $\hskip 0,1 cm \tau:=\tau_{N,P}\in H,\hskip 0,1 cm $  be the transposition defined by \begin{eqnarray*}\hskip 0,1 cm \tau_{N,P}(n)=n, \forall  n\notin\{ N+1,P\},\hskip 1 cm   \tau_{N,P}(N+1)=P\hskip 0,2 cm \textit{ and} \hskip 0,4 cm \tau_{N,P}(P)=N+1. \end{eqnarray*} Then, because $\hskip 0,1 cm \mu\circ T^{-1}_{\tau_{N,P}}=\mu,\hskip 0,1 cm$ for all $\hskip 0,1 cm P\hskip 0,1 cm$  with $\hskip 0,1 cm N+1<P\le Q,\hskip 0,1 cm $ the following equalities 
 \begin{eqnarray*}
\mu(\omega_0\in A_0,...,\omega_{N+1}\in A_{N+1})=\mu(\omega_0\in A_0,...,\omega_N\in A_N, \omega_P\in A_{N+1})\\
={1\over Q-N-1}\sum_{P=N+2}^Q \mu(\omega_0\in A_0,...,\omega_N\in A_N, \omega_P\in A_{N+1})
\end{eqnarray*}
hold  for all $\hskip 0,1 cm Q>N+1,\hskip 0,1 cm$ 
and 
show then, by ergodicity, that
\begin{eqnarray*}
\mu(\omega_0\in A_0,...,\omega_{N+1}\in A_{N+1})= \mu(\omega_0\in A_0,...,\omega_N\in A_N)\mu(\omega_0\in A_{N+1}).
\end{eqnarray*}
The proof is achieved.$\square$\\

  Recall that, if $\mu$ is a probability measure on $\Omega$,  then  a  measurable set  $A$ is $\mu$-almost everywhere  shift-invariant  if $1_A=1_A\circ S$ in $L^1(\mu)$, and also that  it is $\mu$-almost everywhere exchangeable  if $1_A=1_A\circ T_\sigma$ in $L^1(\mu)$, for any $\sigma\in H$.\\Let us use the following notation
  \begin{eqnarray*}
{\cal S}=\textit{ the convex set of all exchangeable probability measures on }\hskip 0,2 cm \Omega=K^\mathbb Z. 
\end{eqnarray*}
   We have seen that $\cal {S}$ is a subset of $M_1(\Omega,S)$, the space of all  shift invariant probability measures on $\Omega$.  It is clear, when  $K$ is a compact space, that $\cal S$ is closed in $M_1(\Omega,S)$ for the weak star topology $\sigma(C(K)^*,C(K))$. Then  ${\cal S}$ is convex compact for this topology. 
We need the following Lemma [15], which  is easy to prove.\\
{\bf{Lemma 3}}\\
\textit{For all $\mu\in{\cal S}$, the sigma-algebra ${\cal I}_\mu$ of $\mu$-invariant sets is equal to the sigma-algebra  ${\cal E}_{xch}$ of $\mu$-exchangeable  sets: 
\begin{eqnarray*}
{\cal I}_\mu={\cal E}_{xch}.
\end{eqnarray*}}
{\bf{Lemma 4}}\\
\textit{ The set of  measures which are extreme points of  ${\cal S}$ is the set of all Bernoulli measures.}\\

{\bf{Proof}} Let $\mu\in {\cal S}$ be  an extreme point. By Proposition 1, it suffices to show that $\mu$ is ergodic for the shift $S$. 
Suppose that $\mu$ is not ergodic. Then there exists an invariant set $A$, with $\mu(A)\in ]0,1[$. Then $A\in {\cal E}_{xch}$ also, because ${\cal I}_\mu={\cal E}_{xch}$, by Lemma 2. Set then
$\mu_1={1\over \mu(A)}1_A\mu$ and $\mu_2={1\over \mu(A^c)}1_{A^c}\mu$. Then $\mu_1, \mu_2\in {\cal S}$, $\mu_1\ne \mu_2$ and $\mu=\mu(A)\mu_1+\mu(A^c)\mu_2$, so that $\mu$ is not extremal in ${\cal S}$, contradicting the hypothesis, and we conclude that $\mu$ is ergodic.$\square$ \\

{\bf{Remark 2}}\\
\textit{It follows from Lemma  4 that the set of exreme points of ${\cal S}$ is closed.}\\

{\bf{Corollary 2}}\\
\textit{ Let $(X_n)$ be an exchangeable sequence with values in a compact metrizable space $K$, with law $\mu$. Then there exists a probability measure $\eta$ supported on the Bernoulli measures such that
for every Borel subset $A$ of $K^\mathbb Z$,
\begin{eqnarray*}
\mu(A)=\int_{Ber}  \tilde{\pi}(A) d\eta(\tilde{\pi})
\end{eqnarray*}
where   $Ber$  stands for the set of all  Bernoulli probability measures on $\Omega=K^\mathbb Z$.
In particular for any $k\ge 1$,  all Borel sets $A_1,...,A_k$ in $K$,
\begin{eqnarray*}
\mu(X_1\in A_1,...,X_k\in A_k)=\int_{Ber} \pi(X_1\in A_1)...\pi(X_1\in A_k)d\eta(\tilde{\pi}).
\end{eqnarray*}}
{\bf{Proof}}\\
  Since $\Omega$ is compact metrizable  $M_1(\Omega,S)$ is compact metrizable. Then  ${\cal S}$ is compact metrizable also, and so, by Choquet Theorem ([18], p.14), there  is a probability measure $\eta$ on ${\cal S}$, which represents $\mu$ and is supported by the extreme points of ${\cal S}$ and this ends the proof by Lemma 4. $\square$\\

{\bf{Remark 3}}\\
\textit{Under the conditions as in the Corollary above, but assuming only the state space $K$ to be compact Haussdorff, the same conclusion holds, with Borel replaced by Baire. }\\ 

{\bf{Corollary 3}}\\
\textit{ If the state space $K$ is compact metrizable, the set  ${\cal S}$ is a simplex.}\\

{\bf{Proof}}
  For any $\hskip 0,1 cm  \mu\in{\cal S},\hskip 0,1 cm $ the decomposition $\hskip 0,2 cm \mu=\int_{{\cal E}} \xi d\eta(\xi),\hskip 0,2 cm $  on the ergodic measures, which is unique since $\hskip 0,2 cm M_1(\Omega,S)\hskip 0,2 cm$ is a simplex, is, by Corollary 2,  given( represented) by  a  measure $\eta$  concentrated on the subset $\hskip 0,2 cm Ber.\hskip 0,2 cm$
Then, because, by Lemma 3,   $Ber$ is the set of extreme points of ${\cal S}$, we conclude that every $\mu\in {\cal S}$ is the barycenter of a unique measure which is supported by the extreme points of $\cal {S} $.\\
Since ${\cal S}$, being a closed subset of the compact metrizable space $M_1(\Omega,S)$, is metrizable,
 it follows then,   by Choquet Theorem(  [ 18], p. 60),  that ${\cal S}$ is a simplex.$\square$\\

\section{Application:
 Bernoullicity of some Determinantal processes}
We recall some properties of simple point processes and determinantal processes. For more details we refer to [5, 19, 20, 22, 4, 10]. 
Let $E$ be a locally compact  Polish space. A locally finite subset of $E$ is called a configuration in $E$. Let $\textit{Conf}\hskip 0,1 cm (E)$ be the space of configurations in $E$. We identify a configuration $\xi$, defined as a set, with the atomic integer valued measure $\sum_{x\in \xi} \delta_x$, on $E$ and then the space $\textit{Conf}\hskip 0,1 cm (E)$ is endowed with the vague topology of measures on $E$, and the corresponding Borel sigma-algebra  ${\cal B}( \textit{Conf}\hskip 0,1 cm (E))$.\\ A simple random point process with phase space $E$  is a probability measure $P$   on the measurable space $(\textit{Conf}\hskip 0,1 cm (E), {\cal B}(\textit{Conf}\hskip 0,1 cm (E))$.  If $B$ is a bounded Borel subset of $E$, let $\eta_B$ be the function defined on $\textit{Conf}\hskip 0,1 cm (E), $ by
\begin{eqnarray*}
\eta_B(\xi):=\textit{card}(B\cap \xi)=\xi(B),\hskip 0,1 cm \forall \xi\in \textit{Conf}\hskip 0,1 cm (E).
\end{eqnarray*}
The joint distribution of the random variables $\eta_B$ determines the measure $P$. 
 Let $\lambda$ be a measure on $E$. A locally integrable function  $\rho^{(n)}$  on the cartesian product $E^n$, is called a $n$-point correlation function of $P$ if
 \begin{eqnarray*}
E_P(\prod _{j=1}^m {\eta_{B(j)}!\over (\eta_{B(j)}-k_j)!})=\int_{B(1)^{k_1}\times ...\times B(m)^{k_m}} \rho^{(n)}(x_1,...,x_n) d\lambda(x_1)...d\lambda(x_n), 
\end{eqnarray*}
for all disjoint bounded Borel subsets  $ B(1),...,B(m)$  of $E$, and all $k_1,...,k_m \in \mathbb N$, with $k_1+...+k_m=n$( [22]  Def. 2 , [4] Def 2, [10])).\\
Another natural way to investigate a point process is to consider expectations of functions $\tilde{\phi}$ of the form
\begin{eqnarray*}
\tilde{\phi}(\xi)=\prod _{x\in \xi} (1+\phi(x))\hskip 2 cm (d_1)
\end{eqnarray*}
 where $\phi$ is a measurable function on $E$, with support in some bounded Borel set. Note that if
  $\mid \phi(x)\mid <1$, for all $x\in E$, then we can write
  \begin{eqnarray*}
\tilde{\phi}(\xi)= \exp(\sum_{x\in \xi} \ln(1+\phi(x))=\exp(\int_E \ln(1+\phi (x))d\xi(x)).
\end{eqnarray*}
Taking, in particular $\phi=-1+\textit{exp}(-\psi)$, with $\psi\ge 0$ having bounded support, we obtain
\begin{eqnarray*}
\tilde{\phi}(\xi)=\exp(\langle -\psi,\xi\rangle)\hskip 1 cm \textit{and then}\hskip 1 cm E_P[\tilde{\phi}]=\int_{\textit{Conf(E)}} \exp(\langle -\psi,\xi\rangle) dP(\xi)
\end{eqnarray*}
leading to Laplace transform of $P$, 
which proves then that $P$ can be characterized by the expectations of $\tilde{\phi}$'s as in $(d_1)$.\\
The point process $P$, all of whose correlation functions $\rho^{(n)}$ exist,  is called  determinantal if there exists a function $  \hskip 0,1 cm k:E\times E\rightarrow \mathbb R$, such that for all $n$ and $x_1,...,x_n\in E$,
\begin{eqnarray*}
\rho^{(n)}(x_1,...,x_n)=\textit{det}[(k(x_i,x_j))_{i,j=1,...,n}].
\end{eqnarray*}
The function $k$ above is called the correlation kernel of the process $P$. It is not unique. For example, if $f(x)\ne 0, \forall x\in E$ and $k'(x,y)={f(x)k(x,y)\over f(y)}$, then 
\begin{eqnarray*}
\textit{det}[(k(x_i,x_j))_{i,j=1,...,n}]=\textit{det}[(k'(x_i,x_j))_{i,j=1,...,n}].
\end{eqnarray*}

In [19, 20] (see also [22] Theorem 3 p. 934) it is proved  that if $K$ is a bounded symmetric integral operator on $L^2(E,\lambda)$ with kernel $k$,  which is also of locally trace class and with spectrum contained in $[0,1]$,  then there exists a unique probability measure $P$ on $\textit{Conf}(E)$, such that for every nonnegative continuous function $\psi$ with compact support, 
\begin{eqnarray*}
\int_{\textit{Conf(E)}} \exp(\langle -\psi,\xi\rangle) dP(\xi)= \textit{Det}(I-\sqrt{(1-e^{-\psi})}K\sqrt{(1-e^{-\psi})})\hskip 2 cm (d_2)
\end{eqnarray*}
  where the determinant is the Fredholm determinant ( see [21]), and where 
the \\operator $\hskip 0,4 cm \sqrt{(1-e^{-\psi})}K\sqrt{(1-e^{-\psi})}\hskip 0,4 cm $ denotes the integral operator with kernel 
\begin{eqnarray*}
L(x,y)=\sqrt{1-e^{-\psi(x)})}k(x,y)\sqrt{1-e^{-\psi(y)})},
\end{eqnarray*}
and moreover the correlation functions of $P$ are given by
\begin{eqnarray*}
\rho^{(n)}(x_1,...,x_n)=\textit{det}(k(x_i,x_j))_{i,j=1,...,n}.
\end{eqnarray*}
In the particular case, where $E$ is countable( in fact $E=\mathbb Z$)  and $\lambda$ is the counting measure, which is relevant to our purpose,  by identifying a subset $A$ of $E$ with its indicator function $1_A$, we can take the configuration space to be equal $\{0,1\}^E $   
and then  $(d_2)$ is equivalent to $(d_3)$ below ([20] Theorem 1.1, see also [12] and [14]  p. 319)
\begin{eqnarray*}
P(\{\omega\in \{0,1\}^E: \omega(e)=1, \forall e\in A\})=P(\{\xi\subset E: A\subset \xi\})=\textit{det}(k(x,y))_{x,y\in A}\hskip 2 cm (d_3)
\end{eqnarray*}
for any finite subset $A$ of $E$.\\
We have the following result\\

{\bf{Theorem 3 }}\\
\textit{Any stationary  discrete time quasi-invariant determinantal process $X=(X_n)_{n\in \mathbb Z}$, with phase space $\mathbb Z$ is isomorphic to a Bernoulli system.}\\

{\bf{Proof}} 
Let $\mu$ be the law of the process $X$. Since $\Omega:=\{0,1\}^\mathbb Z$ is the configuration space of the process  $X$,  $\hskip 0,2 cm \mu$ is a shift invariant probability measure on $\Omega$.  It follows from   Theorem 7 in [22],  that a translation invariant determinantal random point field, with one-particle space $E=\mathbb Z$,  is mixing of any order and then  in particular it is ergodic,  meaning in our setting   that the system  $(\Omega,S, \mu)$ is ergodic.  Hence the process $X$ satisfies the hypothesis in Theorem 1, and therefore  it is isomorphic to a Bernoulli process. $\square$\\
 
In [3] Bufetov considered a class of determinantal processes with phase space $F$, where   $F=\mathbb R$ (the continuous case),  or $F=$  a countable subset of $\mathbb R$, without accumulation points( the discrete case),  corresponding to projection operators with integrable kernels. In the discrete case,  he proved( [3], Theorem 1.6) that they are quasi-invariant, which means  quasi-exchangeable. It follows then  from the preceding Theorem 3, that the   processes in this class, corresponding to the  phase space  $F=\mathbb Z$,  and  which are translation invariant,  are isomorphic to  Bernoulli systems.\\
This applies, in particular, to the discrete sine process. 
Recall that the discrete sine kernel which is translation invariant kernel on the lattice $\mathbb Z$ is defined by
\begin{eqnarray*}
S(x,y;a)=S_1(x-y,a),  \hskip 0,5 cm x,y\in \mathbb Z, \hskip 0,5 cm \textit{where}\\
S_1(x,a):={\textit{sin}( x(\textit{arccos}(a/2))\over \pi x}, \hskip 0,5 cm x\in \mathbb Z, x\ne 0\\
S_1(0,a)={\textit{arccos}(a/2)\over \pi}.
\end{eqnarray*}
where  $a$ is a real number ($-2\le a\le 2$)( [1] p. 486).\\

{\bf{Remark 4}}:\\
\textit{ There exists processes generating  dynamical systems isomorphic to  Bernoulli and which are not quasi-exhangeable}.\\

Before giving some examples, let us, when  $K$ is a finite set $K:=\{0,...,k-1\}$,  denote the cylinder \begin{eqnarray}
C:=\{\omega\in K^\mathbb Z:\omega_0=x_0,...,\omega_n=x_n\}, \hskip 1 cm 
\textit{simply by}\hskip 1 cm  [[x_0,x_1,...,x_n]].
\end{eqnarray}
Then, for any permutation $\hskip 0,1 cm  \sigma \in H, \hskip 0,1 cm $ with $\hskip 0,1 cm  \tau:=\sigma^{-1}$, and such that the support of $\sigma$ is included in $\{0,1,...,n\}$,
\begin{eqnarray*}
T_\sigma^{-1}C=[[x_{\tau(0)},x_{\tau(1)},...,x_{\tau(n)}]]
\end{eqnarray*}
so that, for a stationary  Markov chain defined by a stochastic matrix $\Pi=(\Pi_{i,j})_{i;j=0,...,k-1}$ and an invariant probability vector $p=(p(0),p(1),...,p(k-1))$,
\begin{eqnarray*}
\mu(C)=p(x_0)\Pi_{x_0,x_1}\Pi_{x_1,x_2}...\Pi_{x_{n-1},x_n}\hskip 1 cm  \textit{and}\\
\mu(T_\sigma^{-1}C)=p(x_{\tau(0)}) \Pi_{x_{\tau(0)},x_{\tau(1)}}\Pi_{x_{\tau(1)},x_{\tau(2)}}....\Pi_{x_{\tau(n-1)},x_{\tau(n)}}.
\end{eqnarray*}
(In general,\begin{eqnarray*}
C=\{\omega: \omega_j\in A_j, j\in J\}\Rightarrow T_\sigma^{-1}C=\{\omega:\omega_k\in A_{\tau(k)}, k\in \sigma(J)\}.)
\end{eqnarray*} 
As a simple example,  consider the stationary Markov chain with state space $K:=\{0,1\}$, defined by the matrix  $\Pi$ given by 
 \[
\left(
\begin{array}{cc}
  t&1-t      \\
  1&0      \\
    
\end{array}
\right)
\]

and  the invariant row probability vector $\hskip 0,1 cm  p\hskip 0,1 cm $( meaning  $\hskip 0,1 cm  p\Pi=p$)
\begin{eqnarray*}
p={1\over 2-t}(1,1-t).
\end{eqnarray*}
Now the matrix  $\hskip 0,1 cm \Pi^2\hskip 0,1 cm $ is equal  \[
\left(
\begin{array}{cc}
  t^2+1-t&t(1-t)      \\
  t& 1-t     \\
     
\end{array}
\right)
\]
so that for all $\hskip 0,1 cm  t\hskip 0,1 cm $ with  $\hskip 0,1 cm  0<t<1, \hskip 0,1 cm $  $\hskip 0,1 cm  \Pi\hskip 0,1 cm  $ is irreducible and aperiodic,  and then this Markov chain is isomorphic to a Bernoulli shift.

For $n=2$, for example,  let $\tau=\sigma^{-1}$, be the transposition defined by
\begin{eqnarray*}
\tau(0)=1,\hskip 1 cm \tau(1)=0,\hskip 1 cm \tau(n)=n, \forall n\ne 0,1
\end{eqnarray*}
 and take $C:=[[0,1,1]]$ so that if $D:=T_\sigma^{-1} C$, then $T_\sigma^{-1}D=C$, and hence the following lolds
 \begin{eqnarray*}
\mu(C)\le \mu([[1,1]])=0,\hskip 0,4 cm T_\sigma^{-1}C=[[1,0,1]] \hskip 0,4 cm  \textit{and then}\hskip 0,4 cm  \mu(T_\sigma^{-1}C)\ne 0,\\
\textit{also}\hskip 0,4 cm  \mu(D)\ne 0\hskip 0,4 cm \textit{and}\hskip 0,4 cm  \mu(T_\sigma^{-1}D)=0, 
\end{eqnarray*}
and proves that $\mu$ [resp. $\mu\circ T_\sigma^{-1}$] is not absolutely continuous with respect to $\mu\circ T_\sigma^{-1}$ [resp. $\mu$]

More generally, we have the following\\
{\bf{Proposition 2}}\\
\textit{ Any stationary Markov chain with finite state space $K$, with irreducible and aperiodic transition matrix $\Pi$ having at least one zero entry,  is isomorphic to Bernoulli but it is not quasi-exhangeable.}\\ 

{\bf{Proof}}    Observe first, due to the irreducibility  and aperiodicity of $\Pi$,  that if  $p$ is  the row invariant probability vector, then all the coordinates of $p$ are $>0$.   Let  $i_0,j_0\in K$, such that $\Pi_{i_0,j_0}=0$. 
 Then the  following holds
\begin{eqnarray*}
 \exists a, \Pi_{j_0,a}>0, \hskip 0,2 cm \exists j, \Pi_{i_0,j}>0,\hskip 0,2 cm  \exists n, (\Pi^n)_{a,i_0}>0\\
\textit{that is}\hskip 0,4 cm  \exists x_1,...,x_{n-1}, \hskip 0,3 cm  \Pi_{a,x_1}\Pi_{x_1,x_2}...\Pi_{x_{n-1},i_0}>0.
\end{eqnarray*}
Let, with notation as in (28),  \begin{eqnarray*}
C:=[[j,a,x_1,...,x_{n-1}, i_0, j_0]], \\
D:=[[j_0,a,x_1,...,x_{n-1}, i_0, j]],
\end{eqnarray*}
and $\hskip 0,3 cm \sigma\in  H,\hskip 0,3 cm$ be the transposition  defined by
\begin{eqnarray*}
\sigma(p)=p,\hskip 0,3 cm \forall p\notin\{0,n+2\} \hskip 0,3 cm \textit{and}\\
\sigma(0)=n+2, \hskip 0,5 cm \sigma(n+2)=0. 
\end{eqnarray*}
Then
 \begin{eqnarray*}
T_\sigma^{-1}C=D \hskip 0,3 cm \textit{and} \hskip 0,3 cm T_\sigma^{-1}D=C.
\end{eqnarray*}
But 
\begin{eqnarray*}
\mu(C)=\mu(j)\Pi_{j,a}\Pi_{a,x_1}...\Pi_{x_{n-1},i_0}\Pi_{i_0,j_0}=0, \textit{ and}\\
\mu(D)=\mu(j_0)\Pi_{j_0,a}\Pi_{a,x_1}...\Pi_{x_{n-1},i_0}\Pi_{i_0,j}\ne 0
\end{eqnarray*}
so that, $\mu\circ T_\sigma^{-1}$[ respectively $\mu$]  is not absolutely continuous with respect to $\mu$[ respectively  $\mu\circ T_\sigma^{-1}$], because 
\begin{eqnarray*}
\mu(C)=0,\hskip 0,3 cm \mu(T_\sigma^{-1}C)\ne 0,\\
\mu(T_\sigma^{-1}D)=0, \hskip 0,3 cm \mu(D)\ne 0.
\end{eqnarray*}
The proof is achieved because any mixing Markov chain is isomorphic to a Bernoulli System. $\square$\\
{\bf{Acknowledgments}}.  
Jean Paul Thouvenot posed to me the question  whether or not quasi-invariance implies complete positive entropy. He suggests many improvements and complements. I am greatly indebted to him. I am deeply grateful to him  for many useful discussions as well as for valuable comments which improve  the presentation of the paper,   and also for  bringing to my attention some  appropriate references.
\newpage \centerline{References}
1.  A. Borodin, A. Okounkov and  G. Olshanski(2000): Asymptotics of Plancherel measures for symmetric groups. J. Amer. Math. Soc. 13, 481-515.\\
2.  A. Borodin and  G. Olshanski (2005): Random Partitions and the gamma kernel. Adv. Math. 194,141-202.\\
3.    A. I. Bufetov (2018): Quasi-Symmetries of Determinantal  Point Processes. Ann. Prob., Vol. 46. N$\small {0}$. 2, 956-1003.\\
4.  I. Camilier and L. Decreusefond (2010): Quasi-invariance and integration by parts for determinantal and permanental processes. Journal of Functional Analysis, 259, 268-300.\\
5.  D. J. Daley and   D. Vere-Jones: An Introduction to the Theory of Point Processes, Volume II: General Theory and Structure. Second Edition. Sptinger.\\
6.   N. Dunford and J. T. Schwartz : Linear operators, Part I. General Theory, Interscience publishers, New York.\\
7.  N. A. Friedman and  D. S. Ornstein (1970):  On isomorphism of weak Bernoulli transformations. Adv. in Math. 5, 365-394.\\
8.  E. Hewitt and  L. J. Savage (1955): Symmetric measures on Cartesian products. Trans. Amer. Math. Soc. 80, 470-501.\\
9. E. Hewitt and K. Stromberg: Real and Abstract Analysis. Springer-Verlag, Berlin. Heidelberg. New York(1969).\\
10.  K. Johansson (2005): Random Matrices and Determinantal Processes. arXiv:math-ph/0510038v1 10 Oct 2005.\\
11.  F. Ledrappier (1976):  Sur la Condition de Bernoulli Faible et ses Applications. Lect. Notes in Math., vol 532, 152-159, Springer-Verlag.\\
12.  R. Lyons (2003). Determinantal probability measures. Publ. Math. Inst. Hautes Etudes Sci. 98, 167-212.\\
13.  R. Lyons and  J. E. Steif: Stationary Determinantal Processes: Phase Multiplicity, Bernoullicity, Entropy, and Domination. (2003). Duke Matth. J. Vol. 120, N0. 3, 515-575.\\
14.  G. Olshanski (2008): Difference Operators and Determinantal Point Processes. Functional Analysis and Its Applications, Vol. 42, No. 4,pp.317-329.\\
15.  R. A. Olshen (1971): The Coincidence of Measure Algebras under an Exchangeable Probability. Z. Wahrscheinlichkeitstheorie verw. Geb. 18, 153-158.\\
16.  D. Ornstein  (1970): Bernoulli shifts with the same entropy are isomorphic. Adv. Math. 4,337-352.\\
17. D. Ornstein  (1971):  Two Bernoulli shifts with infinite entropy are isomorphic. Adv. Math. 5, 339-348.\\
18.  R. R. Phelps (2001): Lectures on Choquet's Theorem. Lecture notes in mathematics 1757, 2-ed., Springer, Berlin, New York.\\  
19.  T. Shirai and  Y. Takahashi (2003):  Random point fields associated with certain Fredholm determinants I: fermion, Poisson and boson point processes. Journal of Functional Analysis 205, 414-463.\\ 
20. T. Shirai and  Y. Takahashi (2003):  Random point fields associated with certain Fredholm determinants II: Fermion Shifts and their Ergodic and Gibbs Properties.The Annals of Probability  Vol. 31, N0. 3, 1533-1564. \\
21.  B. Simon (1977): Notes on Infinite Determinants of Hilbert Space operators. Adv. Math. 24, 244-273.\\
22.  A. Soshnikov (2000): Determinantal random point fields. Russ. Math. Surv. 55 923, 107-160.
\\
\\
  {Sorbonne Universit\' e,  UMR 8001, Laboratoire de Probabilit\' es, Statistique et Mod\' elisation,  Bo\^ite courrier 158, 4 Place Jussieu, F-75752 Paris Cedex 05, France.\\
E-mail: doureid.hamdan@upmc.fr}

\end{document}